\documentclass{amsart}
\usepackage{amsthm,graphicx,amssymb}
\usepackage[usenames]{color} 
\usepackage{colortbl}
\usepackage{caption}
\usepackage{subcaption}
\usepackage{hyperref}

\usepackage{lineno}
\usepackage{tikz}
\theoremstyle{definition}
 \newtheorem{theorem}{Theorem}
 
 \newtheorem{definition}{Definition}

\newcommand\RR{\mathbb R}

\newcommand\ZZ{\mathbb Z}
\newcommand\NN{\mathbb N}

\begin{document}

\title{Tropical curves in sandpiles}

\author[N. Kalinin, M. Shkolnikov]{Nikita Kalinin, Mikhail Shkolnikov}

\thanks{
Research is supported in part the grant 159240 of the Swiss National Science
Foundation as well as by the National Center of Competence in Research
SwissMAP of the Swiss National Science Foundation.  
}


\address{Universit\'e de Gen\`eve, Section de
  Math\'ematiques, Route de Drize 7, Villa Battelle, 1227 Carouge, Switzerland}

\email{Nikita.Kalinin\{at\}unige.ch \hfil\break mikhail.shkolnikov\{at\}gmail.com}

\date{\today}

\maketitle

\begin{abstract}
We study a sandpile model on the set of the lattice points in a large lattice polygon. A small perturbation $\psi$ of the maximal stable state $\mu\equiv 3$ is obtained by adding extra grains at several points. It appears, that the result $\psi^\circ$ of the relaxation of $\psi$ coincides with $\mu$ almost everywhere; the set where $\psi^\circ\ne \mu$ is called the deviation locus. The scaling limit of the deviation locus turns out to be a distinguished tropical curve passing through the perturbation points.

Nous consid\'erons le mod\`ele du tas de sable sur l'ensemble des points entiers d'un polygone entier. En ajoutant des grains de sable en certains points, on obtient une perturbation mineure de la configuration stable maximale $\mu\equiv 3$. Le r\'esultat $\psi^\circ$ de la relaxation est presque partout \'egal \`a $\mu$. On appelle lieu de d\'eviation l'ensemble des points o\`u $\psi^\circ\ne \mu$. La limite au sens de la distance de Hausdorff du lieu de d\'eviation est une courbe tropicale sp\'eciale, qui passe par les points de perturbation.

\end{abstract}

\section{Introduction}
\begin{figure}[h]
\minipage{0.32\textwidth}
\fbox{\includegraphics[width=\linewidth]{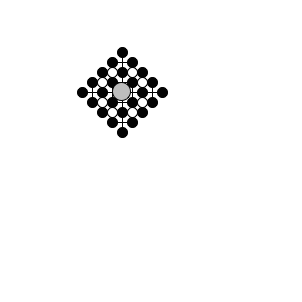}}
\endminipage\hfill
\minipage{0.32\textwidth}
\fbox{\includegraphics[width=\linewidth]{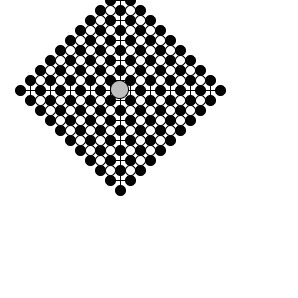}}
\endminipage\hfill
\minipage{0.32\textwidth}
\fbox{\includegraphics[width=\linewidth]{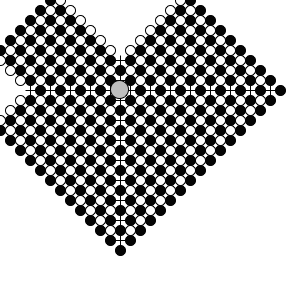}}
\endminipage\hfill

\minipage{0.32\textwidth}
\fbox{\includegraphics[width=\linewidth]{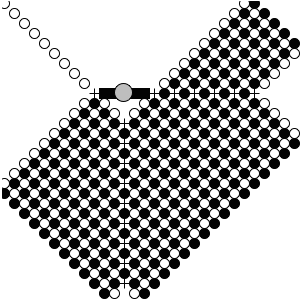}}
\endminipage\hfill
\minipage{0.32\textwidth}
\fbox{\includegraphics[width=\linewidth]{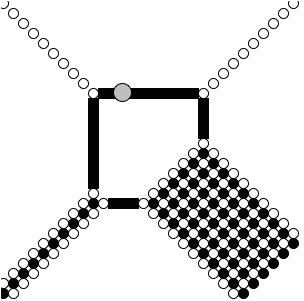}}
\endminipage\hfill
\minipage{0.32\textwidth}
\fbox{\includegraphics[width=\linewidth]{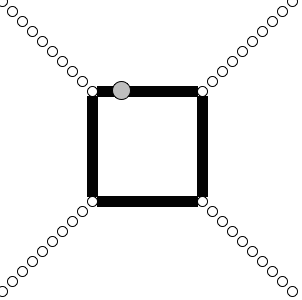}}
\endminipage\hfill
\caption{Snapshots during the relaxation  for the state $\phi \equiv 3$ on a square after adding an extra grain at one point $p$ (the big grey point). Black rounds represent $x$ with $\phi(x)\geq 4$, black squares (which are arranged along the vertical and horizontal edges on the final picture) represent the value of sand equal to $2$, white rounds (arranged along diagonals on the final picture) are $1$, and whites cells are $3$. Rare cells with zero grains are marked as crosses, one can see them during the relaxation on the vertical and horizontal lines through $p$. The value of the final state at $p$ is $3$.\medskip \\
Instantan\'es pendant la relaxation de la configuration $\phi\equiv 3$ sur un carr\'e, apr\`es ajout d'un grain additionnel au point p (le gros point gris).
Les ronds noirs repr\'esentent x avec $\phi(x)\geq 4$, les carr\'es noirs (qui se trouvent sur les ar\^etes verticales et horizontales dans la configuration finale)
repr\'esentent les points o\'u le nombre de grains de sable \'egale $2$. Les cercles (sur les diagonales dans la configuration finale) poss\`edent un grain de sable et les cellules blanches en ont $3$. Les cellules rares avec z\'ero grains sont marqu\'ees par des croix; on peut les voir pendant la relaxation sur les lignes verticales et horizontales passant par $p$. La valeur au point p dans la configuration finale est $3$.} 
    \label{fig_onepoint}
\end{figure}

\subsection{Sandpile model on a finite set} 
 Consider the standard lattice $\ZZ^2$ on the plane. For $v\in{\ZZ^2}$ we denote by $n(v)$ the set of all four closest points to $v$ in $\ZZ^2.$ Let $\Gamma$ be a finite subset of $\ZZ^2.$ A {\it state} (or a {\it configuration}) $\phi\in\ZZ^\Gamma$  of a sandpile on $\Gamma$ is a non-negative integer-valued function on $\Gamma.$ For a state $\phi$ and $v\in\Gamma$ we interpret $\phi(v)$ as a number of sand grains at the site $v.$ For each $v\in\Gamma$ we define the {\it toppling} operator $T_v\colon \ZZ^\Gamma\rightarrow\ZZ^\Gamma$ at $v$ given by $$T_v\phi=\phi-4\delta_v+\sum_{w\in \Gamma\cap n(v)}\delta_w,$$ where $\delta_v$ is a function on the lattice defined to be $1$ at $v$ and $0$ otherwise.  A toppling $T_v$ is called {\it legal} for a state $\phi$ if $\phi(v)\geq 4$, i.e. if $T_v\phi$ is again a state.  In this case we think of $T_v$ as a redistribution of sand from the overfilled site $v$ to its neighbors. If some neighbors are missing in $\Gamma$, i.e. $n(v)\not\subset\Gamma$, then at least one grain leaves the system after the toppling.  If $\phi(v)<4$ for all $v\in\Gamma$, then $\phi$ is called a {\it stable} state. The state $\mu$ which is defined to be equal to $3$ at every point of $\Gamma$ is called the {\it maximal stable state}.

A {\it relaxation} for a state $\phi$ is a sequence of states $\phi=\phi_0,\phi_1,\dots\phi_m$ such that $\phi_{i+1}$ is the result of applying a legal toppling to $\phi_i$ and $\phi_m$ is a stable state. It is well known that for any state $\phi$ there exists a relaxation and the last state $\phi_m$ depends only on $\phi$ (see \cite{BTW,Dhar}). We denote $\phi_m$ by $\phi^\circ$ and call it {\it the result of the relaxation} of $\phi.$ Informally, in order to find the result of the relaxation it doesn't matter which legal topplings to apply at each step of a relaxation sequence.  
\subsection{Motivation}
Let $\Omega\subset\RR^2$ be a non-degenerate (i.e. of non-zero area) lattice polygon.  Let $\Gamma$ be the intersection of $\ZZ^2$ with $\Omega.$ Consider a set $P\subset \Gamma$. We add extra grains to the state $\mu$ at the points $P=\{p_1,p_2,\dots, p_n\}$. After the relaxation, this gives the state $\phi=(3+\sum_{p\in P}\delta_p)^\circ.$  

 In the examples shown in Figures \ref{fig_onepoint},\ref{fig_triang} and \ref{fig_square} we see that the set of points where $\phi$ is not maximal constitutes some sort of a graph passing through $P$. As it is seen in Figures \ref{fig_triang} and \ref{fig_square}, the picture becomes more regular when the cardinality of $P$ is small with respect to the size of $\Omega.$  In the next section we state certain theorems formalizing this concept. Some of this results and their far going generalizations will appear in \cite{us}. This short note can be seen as an introduction to the subject and an announcement of these results.

\subsection{Related activities}What we described above is quite close, at least at the basic level, to the original Back-Tang-Wiesenfeld model \cite{BTW}. The fundamental difference is that their framework had a probabilistic flavor. Their subject was a stochastic growth of a sandpile simulated by an iterative process of adding an extra grain at a random site followed by a relaxation. 

In the paper \cite{Dhar} it has been proven that for any initial configuration most of the states stop occurring after some step of the Back-Tang-Wiesenfeld-process and the other states (the {\it recurrent} states) occur with equal probabilities. Dhar has also shown that the set of all recurrent states forms an abelian group (called the sandpile group) under the action $\phi\oplus\psi=(\phi+\psi)^\circ.$ Study of this group is an important and fruitful part of the theory of abelian sandpiles \cite{LP}. Unfortunately, it is little known about the structure of particular elements of the sandpile group \cite{BR}. There is a simple characterization of a recurrent state. A stable state $\phi$ is recurrent if and only if it can be represented as $\phi=(3+\psi)^\circ$ for some state $\psi.$ Therefore, the states that we study in this paper are recurrent.

Our initial inspiration was the definite presence of a tropical curve in Figure 2 (which coincide with our Figure~\ref{fig_onepoint}) in \cite{firstsand}, where the authors considered the case where $\Omega$ is a square and $P$ consists of just one point.  The same authors experimentally found more pictures in Section 4.3 of \cite{CPS} (Figure 3.1 in \cite{book}) for the case of others $\Omega$ and $P$; these pictures contain tropical curves too. Note that they considered a bit different process: after dropping a grain of sand to a point $p$ and subsequent relaxation, they remove a grain from $p$.  After a number of computer simulations, it appeared that the toppling function for the relaxation of $3+\sum_{p\in P}\delta_p$ is almost piece-wise linear, it was also observed in \cite{sadhu2012pattern}, where it was conjectured that tropical geometry could be useful in this problem. These articles, \cite{firstsand,CPS}, containe a lot of observations and heuristic arguments. In order to formulate these observations rigorously, we use scaling limits. It makes our results particularly close to \cite{FLP, LPS, PS} where the scaling limit of the states $(N\delta_0)^\circ$ was shown to exist. Then, the proofs appeared to be rather cumbersome, but evolving into a deep connections between sublinear order integer-valued discrete harmonic functions, dynamics of polytopes, and tropical geometry, see \cite{us} for details.

\section{The results}

\begin{figure}
    \centering
    \setlength\fboxsep{0pt}
    \setlength\fboxrule{0.5pt}
    \includegraphics[width=0.3\textwidth]{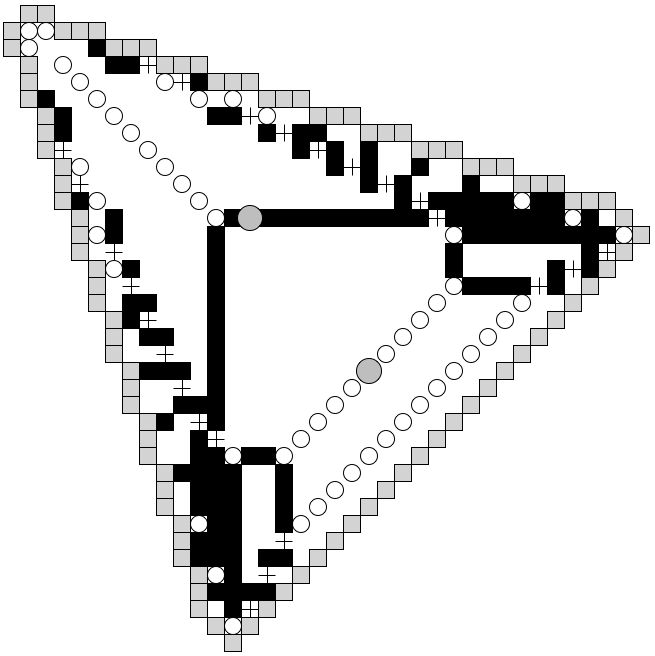}\ 
    \includegraphics[width=0.3\textwidth]{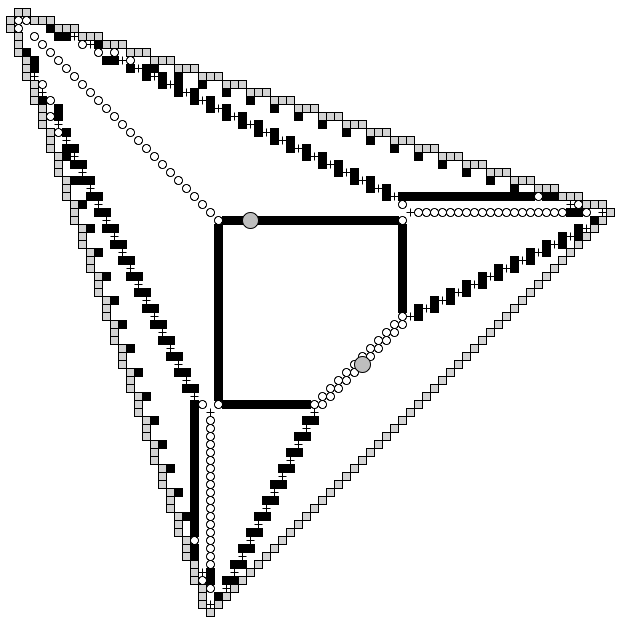}\ 
    \includegraphics[width=0.3\textwidth]{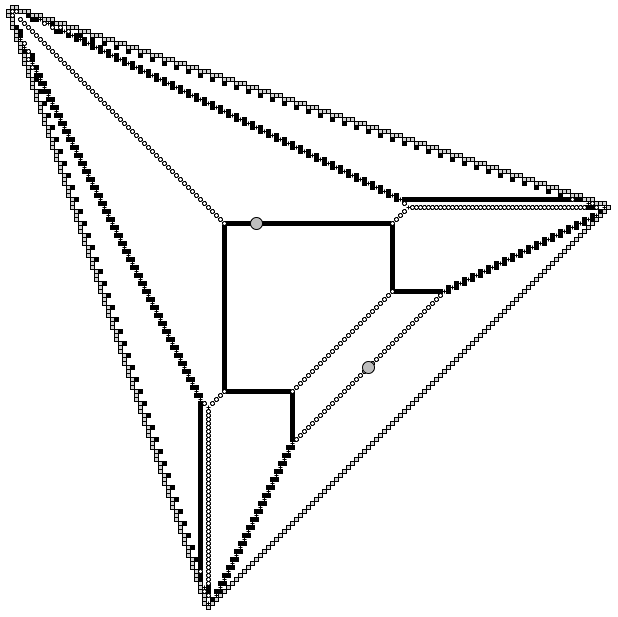}\  
    \caption{Triangular domain and two perturbation points. The scale of last picture is four times bigger than the scale of the first and two times bigger than the scale of the second.\\  Un domaine triangulaire, deux point de perturbation. L'\'echelle de la figure de droite est quatre fois celle de la figure de gauche et deux fois celle de la figure centrale. } 
    
    \label{fig_triang}
\end{figure}

 Let $\Omega$ be a non-degenerate lattice polygon and $P$ be a finite non-empty subset of $\Omega^\circ.$ For any $N\in\NN$ consider a set $\Gamma_N=N\Omega\cap\ZZ^2.$ Denote by $[p]\in\ZZ^2$ the coordinate-wise rounding down of a point $p\in\RR^2$. Define the state $\phi_N=(3+\sum_{p\in P}\delta_{[Np]})^\circ$ on $\Gamma_N$ and the {\it deviation} set $$C_N={1\over N}\{v\in\Gamma_N | \phi_N(v)<3\}.$$  Experimental evidence suggests that when $N$ grows, the shape of sets $C_N\subset \Omega$ stabilizes, see Figure \ref{fig_triang}. 
 
\begin{theorem}\label{th_polygoneconvergence}
 The sequence of sets $C_N\subset\Omega$ Hausdorff  converges to a set $\tilde C_{\Omega,P}.$ The set $\tilde C_{\Omega,P}$ is a planar graph passing through the points $P$. Each edge of $\tilde C_{\Omega,P}$ is a straight segment with a rational slope. \end{theorem}

Denote by $C_{\Omega,P}$ the closure of $\tilde C_{\Omega,P}\cap\Omega^\circ.$ We have a canonical mapping associating the graph $C_{\Omega,P}$ to the configuration of points $P$. Theorem \ref{th_mimpoly} later represents $C_{\Omega,P}$ as a solution of some variation problem.

First of all we note that $C_{\Omega,P}$ is a weighted graph, i.e. there is a canonical choice of weights for its edges. This choice comes from averaging an amount of sand in $\phi_N$ along edges of $C_{\Omega,P.}$  Namely, we define a sequence of functions $\phi'_N\colon\RR^2\rightarrow\mathbb{R}$ given by $\phi'_N(x)=N(3-\phi_N([Nx]))$ if $[Nx]\in\Gamma_N$ and $\phi'_N(x)=0$ otherwise.
 
 \begin{theorem}\label{th_limphi}
 There exists a $*$-weak limit $\tilde\phi$  of the sequence $\phi'_N.$ Moreover, there exists a unique assignment of weights $m^{\Omega,P}_e\in\NN$ for the edges of $C_{\Omega,P}$ such that for all smooth functions $\psi$ supported on $\Omega^\circ$ $$\tilde \phi(\psi)=\lim_{N\to\infty}\int_{\RR^2} \phi'_N \psi=\sum_{e\in E}||l_e||m^{\Omega,P}_e\int_e \psi,$$ where $E$ is the set of all edges of $C_{\Omega,P}$ and  $l_e$ is a primitive vector of $e\in E$, i.e. the coordinates of $l_e$ are coprime integers and $l_e$ is parallel to $e.$ 
\end{theorem}

\begin{figure}
    \centering
    \setlength\fboxsep{0pt}
    \setlength\fboxrule{0.5pt}
    
    \includegraphics[width=0.45\textwidth]{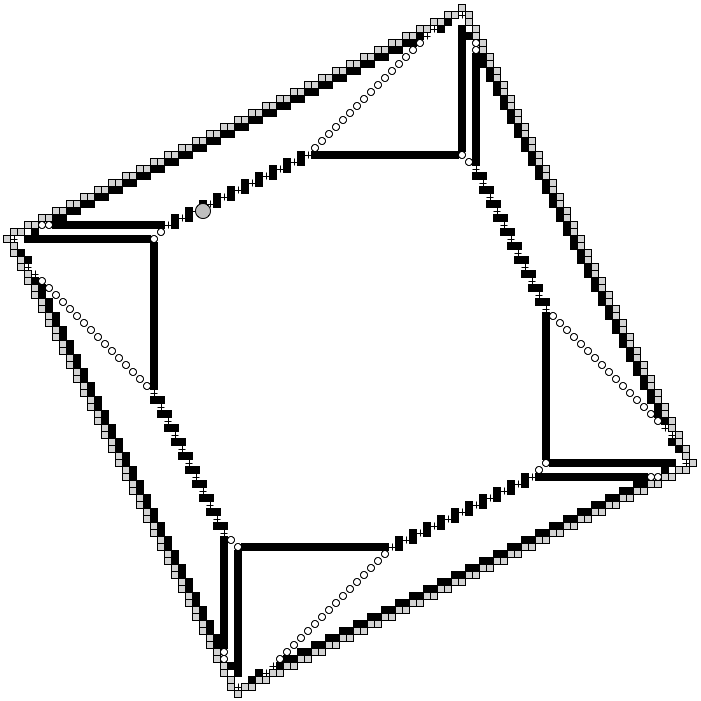}
    \begin{picture}(200,100)
   
\put(0,0){\includegraphics[width=0.45\textwidth]{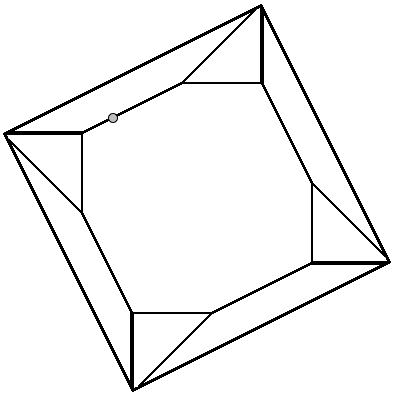}\ }
\put(61,30){$2$}
\put(29,140){$2$}
\put(173,73){$2$}
\put(139,169){$2$}
\put(139,169){$2$}
\end{picture}\hfil
    \caption{The result of adding a grain in two different scales. The polygon $\Omega$ is a square with edges having slopes $2$ and $1/2.$ \\ R\'esultat de l'ajout d'un grain de sable \`a deux \'echelles diff\'erentes. Le domaine $\Omega$ est un carr\'e avec c\^ot\'es de pente $2$ et $1/2$ respectivement.} 
    \label{fig_square}
\end{figure}

In other words, edges with high weights are ``wider'' than edges with low weights (see Figures \ref{fig_triang} and \ref{fig_square}).

Theorem \ref{th_limphi} motivates the following definition.  

\begin{definition}[\cite{tony}]
Let $C$ be a weighted planar graph whose edges have rational slopes. We define {\it the tropical symplectic area} of $C$ by $$Area(C)=\sum_{e\in E} ||l_e|| m_e ||e||.$$
\end{definition}

It appears that $C_{\Omega,P}$ minimizes its area in the class of weighted graphs that we call $\Omega$-{\it tropical} curves. 

\begin{definition}
A finite planar weighted graph $C\subset\Omega$ is called an $\Omega$-tropical curve  if \begin{itemize}
\item each edge $e$ of $C$ has a rational slope and  an integer weight $m_e,$ 
 \item the intersection of $C$ with $\partial\Omega$ is the set of all vertices of $\Omega$,
 \item if $v\in\Omega^\circ$ is a vertex of $C,$ then \begin{equation} \label{eq_balancing}\sum_{e\in E_v} m_e l_e=0, \end{equation}  where $E_v$ is the set of all edges of $C_{\Omega,P}$ incident to $v$ and $l_e$ is the primitive vector of $e\in E_v$ oriented out of $v,$
 \item there exist a unique labeling $d_s\in\ZZ_{>0}$ for each side $s$ of $\Omega$ such that for each vertex $v$ of $\Omega$ \begin{equation}\label{eq_outerbalancing}\sum_{e\in E_v} m_e l_e=d_{s_1}l_{s_1}+d_{s_2}l_{s_2}, \end{equation}
where $s_1$ and $s_2$ are the two sides of $\Omega$ incident to $v,$ $l_{s_i}$ is the primitive vector of $s_i$ oriented out of $v.$
 \end{itemize}
 \end{definition}

Condition \eqref{eq_balancing} is well known in tropical geometry under the name {\it balancing condition,} see \cite{mikh}. We call \eqref{eq_outerbalancing} the {\it outer balancing} condition.

The condition \eqref{eq_balancing} implies that an $\Omega$-tropical curve can be always represented as the intersection of $\Omega$ with a {\it tropical curve} (see \cite{mikh,BIMS}).  A {\it plane tropical curve} \label{tropcurvedef} is a planar graph, whose edges are intervals with rational slopes and prescribed positive integer weights. The edges are allowed to be unbounded, the number of edges and vertices is finite and at each vertex the balancing condition is satisfied. If $C$ is a subgraph of a tropical curve, then its tropical symplectic area can be seen as the limit of the usual symplectic areas for a family of holomorphic curves that degenerates to $C$ (see \cite{us,tony}). On the other hand, this area can be seen as a particular normalization for the Euclidean length of $C.$ In this interpretation $C_{P,\Omega}$ gives a solution to the analog of the Steiner tree problem.

\begin{theorem}\label{th_minarea}
The graph $C_{\Omega,P}$ is an $\Omega$-tropical curve. Furtherefore, $C_{\Omega,P}$ has the minimal tropical symplectic area among all $\Omega$-tropical curves passing through the configuration of points $P.$
\end{theorem}

In certain cases, the curve $C_{\Omega,P}$ can be found by means of this minimization property. This happens essentially when $\Omega$ is ``simple" enough and the number of points $|P|$ is small. In general, for any lattice polygon $\Omega$ there are {\it generic} configurations $P$ for which $C_{\Omega,P}$ is not a unique minimizer of the symplectic area. Fortunately, the curve $C_{\Omega,P}$ can be characterized by the property of minimizing its tropical polynomial.

A {\it tropical polynomial} is a function $F\colon \RR^2\rightarrow \RR$ given by $$F(x_1,x_2)=\min_{(k,l)\in A} (kx_1+lx_2+a_{k,l}),$$ where $A$ is a finite subset of $\ZZ^2$ and $a_{k,l}\in \RR.$ For a tropical polynomial $F$ we consider its corner locus $C,$ i.e. the set of points where $F$ is not locally linear. The set $C$ has a natural structure of a tropical curve \cite{mikh} and we say that the curve $C$ is defined by the polynomial $F$. 

Consider an $\Omega$-tropical curve $C.$ Consider an extension of $C$ to a tropical curve $C',$ i.e.  $C=C'\cap\Omega.$ The outer balancing condition \eqref{eq_outerbalancing} implies that there exists a unique tropical polynomial $F'$ vanishing at $\partial\Omega$ which defines $C'.$ It is easy to see that $F=F'|_\Omega$ does not depend on $C'.$ Therefore, we call $F$ {\it the $\Omega$-tropical polynomial of the curve} $C.$ Denote by $F_{\Omega,P}$ the $\Omega$-tropical polynomial of $C_{\Omega,P}.$ 

\begin{theorem}\label{th_mimpoly}The polynomial $F_{\Omega,P}$ is the point-wise minimum of all polynomials $F$ of all $\Omega$-tropical curves that pass through $P.$ 
\end{theorem}

Theorem \ref{th_mimpoly} can be seen as a manifestation of the least action principle for sandpiles \cite{FLP}. Consider a state $\phi$ on a finite set $\Gamma\subset\ZZ^2.$  Consider a set $\mathcal{H}$ of functions $H\colon\ZZ^2\rightarrow\ZZ_{\geq 0}$ such that $H$ vanishes outside $\Gamma$ and $\phi+\Delta H\leq 3$ on $\Gamma,$ where $\Delta H$ is the discrete Laplacian of $H$ given by $$\Delta H(v)=-4H(v)+\sum_{w\in n(v)} H(w).$$  Denote by $F$ the point-wise minimum of all functions $H\in\mathcal{H}.$ Then $\phi^\circ=\phi+\Delta F.$ Furthermore, $F(v)$ counts the number of topplings at $v\in\ZZ^2$ in any relaxation sequence for $\phi.$  Therefore, we call the function $F$ the {\it toppling function} (or odometer) of the state $\phi.$

Let $F_N$ be the toppling function  for the relaxation of $3+\sum_{p\in P}\delta_{[Np]}$ on $\Gamma_N=N\Omega\cap\ZZ^2.$ Consider a function $\tilde F_N\colon\Omega\rightarrow\RR$ given by $$\tilde F_N(x)={1\over N}F_N([Nx]).$$

\begin{theorem}\label{th_topplingconv}
The tropical polynomial $F_{\Omega,P}$ is the scaling limit of toppling functions, i.e. the sequence $\tilde F_n$ converges pointwise to $F_{\Omega,P}.$ 
\end{theorem}

This theorem essentially implies all the previous ones.  The main idea of the proof of Theorem \ref{th_topplingconv} is that the functions $\tilde F_N$ are bounded by the piecewise linear function $F_{\Omega,P}$ and are harmonic on the regions where $\phi_N=\mu$. Then, $F_N$ is harmonic almost everywhere for large $N$, and its growth is tied by linear functions. This finally implies that $\tilde F_N$ is a piecewise linear function almost everywhere. Detailed proofs and generalizations of these results will appear in the paper \cite{us}. In this forthcoming paper, the lattice polygon $\Omega$ will be replaced by an arbitrary convex domain. In this case, all the convergence results still hold without significant changes. The set $C_{\Omega,P}$ is still a union of straight edges, but the number of edges is usually infinite. It appears that $C_{\Omega,P}$ is a {\it tropical analytic curve} and $F_{\Omega,P}$ is a {\it tropical series}.

\bibliography{../sandbib}

\bibliographystyle{abbrv}
\end{document}